\documentclass[12pt]{article}
\usepackage{latexsym}
\usepackage{amsmath,amsfonts,amssymb,epsfig}
\usepackage{amsthm}
\usepackage{hyperref}

\newtheorem{theorem}{Theorem}[section]
\newtheorem{definition}[theorem]{Definition}
\newtheorem{corollary}[theorem]{Corollary}

\newtheorem{lemma}[theorem]{Lemma}
\newtheorem{conjecture}[theorem]{Conjecture}
\newtheorem{claim}[theorem]{Claim}
\newtheorem{rem}[theorem]{Remark}
\newenvironment{remark}{\begin{rem} \rm}{\end{rem}}
\newtheorem{ide}[theorem]{Idea}

\newtheorem{exa}[theorem]{Example}

\newtheorem{ques}[theorem]{Question}

\newtheorem{pro}[theorem]{Problem}
\newenvironment{problem}{\begin{pro} \rm}{\end{pro}}
\newtheorem{spe}[theorem]{Speculation}

\newcommand{\Lie}{{\cal L}}

\font\bbb=msbm10 scaled 1100
\newcommand{\real}{\mbox{\bbb R}}       

\newcommand{\bfx}{{\mbox{\bf x}}}

\title{On the contact geometry of nodal sets}
\bigskip
\bigskip
\author{R. Komendarczyk\footnote{School of Mathematics,
Georgia Institute of Technology, Atlanta GA, 30332; e-mail: {\it
rako@math.gatech.edu},  research partially supported by  {\it NSF} grant
{\it DMS-0134408}}}

\begin{document}

\maketitle

\begin{abstract}
\quad In the 3-dimensional Riemannian geometry, contact structures equipped with
an adapted Riemannian metric are divergence-free, nondegenerate
eigenforms of the Laplace-Beltrami operator.
We trace out a 2-d analogue of this fact:
there is a close relationship between the topology of the contact structure
on a convex surface in the 3-manifold (the {\em dividing curves}) and
the nodal curves of Laplacian eigenfunctions on that surface.
Motivated by this relationship, we consider a topological version of Payne's
conjecture for the free membrane problem.
We construct counterexamples to Payne's conjecture for closed Riemannian
surfaces.
In light of the correspondence between the nodal lines and
dividing curves, we interpret Payne's conjecture
in terms of the  {\it tight} versus {\it overtwisted}
dichotomy for contact structures.
\end{abstract}
\bigskip
\begin{center}
\small {\bf keywords}: nodal lines, dividing curves, contact structures,
eigenfunctions of Laplacian.
\end{center}

\section{Introduction.}
\subsection{Payne's conjecture.}
\qquad If we think of a given Riemannian surface $(\Sigma, g_\Sigma)$ as a
vibrating membrane with $u(\bfx,t)$ a displacement of the membrane from the
original position in time $t$, $u$ is a solution to the wave equation
\begin{gather}\label{eq:wave}
  \partial_{tt} u=\Delta_{\Sigma} u,
\end{gather}
Since the solution is separable, i.e. $u(t,\bfx)=v(t)w(\bfx)$,
we obtain an equivalent system of
equations $\partial_{tt} v= \lambda v$ and $\Delta_\Sigma w= \lambda w$,
($\lambda\in \real$). Therefore, the ``stagnation points'' on the membrane
are exactly zeros of the eigenfunction $w$. This zero set,
$\Xi(w):=\{\bfx\in\Sigma: w(\bfx)=0\}$, is called a {\it
nodal set} and forms interesting patterns, as originally studied by E.
Chladni in the 18th century.  In case when the membrane is a closed surface
we refer to the problem \eqref{eq:wave} as the {\it free membrane problem},
for surfaces with boundary and Dirichlet boundary conditions we refer to
the problem as the {\it fixed membrane problem} (see \cite{Cheng76}.) For
an arbitrary smooth Riemannian surface $(\Sigma,g_\Sigma)$, the structure
of nodal sets has been characterized by S. Cheng in \cite{Cheng76}, where
it is proved that the nodal set is a collection of $C^2$- immersed closed
curves in $\Sigma$. For a generic metric, K. Uhlenbeck
\cite{Uhlenbeck76} showed that these curves are embedded circles with no
critical points. Not much is known about the general structure of such
sets. One of the fundamental results is Courant's nodal domain theorem. It
states that the nodal set of the $k$th eigenfunction of the Laplacian
divides a domain into at most $k$ regions in case of the fixed membrane
problem, and $k+1$ regions in case of the free membrane problem. Courant's
theorem implies that the first (second) eigenfunction of $\Delta_\Sigma$
has to divide the region into exactly two domains for the
free (fixed) membrane problem. In \cite{Payne67}, L. E. Payne conjectured
that in case of the fixed membrane problem for bounded domains in
$\real^2$, the second eigenfunction of Laplacian cannot posses a closed
nodal line.

\begin{conjecture}[Payne(1967)]
\label{co:payne}
  The second eigenfunction of the Laplacian on a bounded region
  $\Omega$ in Euclidean $\real^2$ with
  the Dirichlet boundary conditions cannot have a closed loop in
  its nodal set.
\end{conjecture}

Since 1967, Payne's conjecture has been proved to be true in the case of
convex domains (see \cite{Alessandrini94} and \cite{Melas92}). Recently, it
has been proved false by T. Hoffmann-Ostenhof and co-authors (see
\cite{Hoffmann-Ostenhof97}), in the case of a non-simply connected domain
(disk with slits on an inner circle removed). It is still not known,
however, if Conjecture \ref{co:payne} is true for an arbitrary simply
connected region in $\real^2$. In \cite{Freitas02}, P. Freitas has shown
that Conjecture fails in case of $\Omega=D^2$ for a non-Euclidean metric.

We consider a more global version of Payne's conjecture:
\begin{problem}\label{co:payne_closed}
  Does the first $\Delta_\Sigma$- eigenfunction on a given closed surface
    $\Sigma$ admit a contractible nodal curve in $\Sigma$?
\end{problem}
The principal result of the paper is the construction of examples
which answer this question in the affirmative.

We observe that Problem~\ref{co:payne_closed} is closely related
to ideas coming from the topology of {\em contact structures}: in
dimension three, these are fields of 2-d planes which are maximally
nonintegrable --- they are as far away from defining a foliation
as possible. The topology of contact structures is extremely
interesting and has been the focus of much research among topologists.
In particular, much progress has been made through the elucidation
of a dichotomy between the {\em tight} and {\em overtwisted}
contact structures (see Section~\ref{sec:Contact structures and their dividing
curves} for definitions).

Our principal observation is a connection between certain topological
features of a contact structure --- so-called {\em dividing curves}
associated to a {\em convex surface} $\Sigma$ --- and the nodal
sets for eigenfunctions of the Laplacian on $\Sigma$. These
results yield a reformulation of Problem~\ref{co:payne_closed}:
\begin{problem}
\label{co:payne_closed_2}
Given $(\Sigma,g)$, is the contact structure on $\Sigma\times\real$
induced by the first eigenfunction of $\Delta_\Sigma$ overtwisted?
\end{problem}

These two conjectures are cousins and provide some basis for
a spectral geometry interpretation of the tight-overtwisted
dichotomy in contact topology.

\subsection{Outline and terminology.}
In Section
\ref{sec:Contact structures and their dividing curves}, we give a short
overview of contact topology and introduce the relevant definitions.
The next section is devoted to the relationship between nodal sets
and dividing curves in contact geometry.
Section \ref{sec:Payne's conjecture for closed Riemannian  surfaces}
carefully constructs counterexamples to Payne's conjecture on closed
Riemannian surfaces. Specifically, we prove that an arbitrary
orientable surface admits metrics such that the principal eigenfunction
has nodal set a single closed
contractible curve. Our technique is based on the work of J. Takahashi,
\cite{Takahashi02}, about collapsing connected sums of surfaces, which is
in turn based on work of C. Anne (see \cite{Anne93}). As an additional
result, we show the $C^\infty$- convergence of eigenfunctions on compact
subsets of the ``non-collapsing'' part of the manifold.

Here, all manifolds, unless stated otherwise, are equipped with a
Riemannian metric, and are compact smooth orientable of dimension either
two or three. Throughout the article $C^j(M)$ stands for the set of $j$-
differentiable functions on $M$, with $j=\infty$ smooth, and $j=\omega$
analytic. Spaces $L^2(M)$, $H^j(M)$ are customary, square integrable real
functions, and the Sobolev space of real valued functions with at least $j$
bounded weak derivatives. The space $\Omega^k(M)=C^\infty(\Lambda^k M)$ is
a set of smooth real valued $k$- differential forms on $M$ making
$\Omega^\ast(M)=\bigoplus^n_{k=0} \Omega^k(M)$ a graded $C^\infty(M)$
module over $\real$, where $n=\dim(M)$. Here we denote by $L^2(\Lambda^k
M)$ and $H^j(\Lambda^k M)$, respectively, the square integrable, and the
Sobolev spaces of $k$- differential forms, where the measure is induced
from the Riemannian metric. The Riemannian metric also induces an $L^2$-
isometry: $\ast:\Omega^k(M)\to\Omega^{n-k}(M)$, namely the {\it Hodge star}
operator. Consequently, we obtain de'Rham graded complexes
$(\Omega^\ast(M), d)$ and  $(\Omega^\ast(M), \delta)$,
where $d\equiv d^k:\Omega^k(M)\to\Omega^{k+1}(M)$ is an exterior derivative
(also called a \emph{differential}), and
$\delta\equiv\delta^k:\Omega^{k+1}(M)\to\Omega^k(M)$ an adjoint of $d$
(also called a \emph{co-differential}) given in terms of the Hodge star by
$\delta^k=(-1)^{k+1}\ast d^{n-k-1}\ast$ or equivalently as a formal adjoint
of $d$,
\begin{gather*}
  (d^k\omega,\eta)_{L^2(\Lambda^{k+1}
M)}=(\omega,\delta^k\eta)_{L^2(\Lambda^k M)},\qquad \omega\in\Omega^k(M),
\eta\in \Omega^{k+1}(M).
\end{gather*}
Most of the time we skip the superscripts in the notation for
differentials and co-differentials and simply write $d$ and $\delta$. The
Laplacian on $k$- forms is defined by $\Delta=\delta\,d+d\,\delta$, which
in the case of functions reduces to $\Delta=\delta\,d$ (for further
reference consult \cite{Rosenberg97} or \cite{Aubin98}).

We also introduce the following notation for nodal sets. Let $\Xi(M,f)=\{x\in
M: f(x)=0\}$ stand for the zero set of the function $f$. In the case
$f=f_k$,  where $f_k$ is $k$th- eigenfunction of $\Delta_M$, we write
$\Xi(M,k):=\Xi(M,f_k)$, or $\Xi(M):=\Xi(M,1):=\Xi(M,f_1)$ for $k=1$.

\section{Contact structures and their dividing curves.}
\label{sec:Contact structures and their dividing curves}
\qquad
Let $M$ be a three dimensional, closed oriented manifold.
A smooth plane field $\xi$ on $M$ is called a {\em contact structure}
if $\xi$ is maximally nonintegrable;
that is, for any pair of vector fields $X$ and $Y$ satisfying $X_p, Y_p\in
\xi_p\subset T_p M$ locally,  we have $[X,Y]_p\notin \xi_p$. (This
condition is exactly opposite to the Frobenius condition for  integrable
subbundles.) The plane field $\xi$ can always be defined as the kernel of
differential 1-form $\alpha$ on $M$. The 1-form $\alpha$ is determined up to
a multiplication by a positive real valued function and is called a {\it
contact form}.
In terms of $\alpha$ the condition of non-integrability can be
expressed as follows,
  \begin{gather}\label{eq:contact_condition}
   \alpha\wedge d\alpha\neq 0.
  \end{gather}
It was proved some time ago by R. Lutz and J. Martinet (in
\cite{Martinet71}) that every closed 3-manifold admits a contact
structure. Since then, there has been a significant amount of research
devoted to the problem of classifying contact structures up to an isotopy
of plane fields (see e.g. \cite{Etnyre00}, \cite{Honda00},
\cite{Honda00_1}). One of the fundamental results in this direction is a
theorem of Y. Eliashberg \cite{Eliashberg89} which divides contact structures
into two classes: {\it overtwisted} and {\it tight}.
\begin{definition}
A contact structure $\xi$ is overtwisted if and only if there
exists an embedded disk $D^2\subset M$ such that $D$ is transverse to $\xi$
near $\partial D$ but $\partial D$ is tangent to $\xi$.
Any contact structure which is not overtwisted is called tight.
\end{definition}
In \cite{Eliashberg89}, Y. Eliashberg classified overtwisted contact
structures in terms of the homotopy type of plane fields. On the other
hand, the complete classification of tight structures still remains an open
problem. In the study of this problem the concept of {\it dividing curves}
for {\it convex surfaces} plays a major role
\cite{Kazez02}, \cite{Honda00}, \cite{Honda00_1}.

\begin{definition}
A {\it convex surface} is a properly embedded surface $\Sigma$ in $(M,\xi)$
such that there exists a vector field $v\pitchfork \Sigma$ transverse to
$\Sigma$ and preserving $\xi$ (i.e. $\Lie_v\xi=0$). The vector field
$v$ is called a {\it contact vector field}.
The {\it dividing set}, $\Gamma_\Sigma$, is the set of all points $p$ on
the surface $\Sigma$ where $v_p\in\xi_p$.
\end{definition}

The following theorem by Giroux (in \cite{Giroux91})
characterizes the dividing set $\Gamma_\Sigma$.
\begin{theorem}
Let $\Sigma$ be a convex surface in $(M,\xi)$. The dividing set
$\Gamma_\Sigma$ of $\xi$ is a set of smooth curves. Moreover, the isotopy
type of $\Gamma_\Sigma$ is independent of choice of the contact field $v$.
\end{theorem}
A parallel theorem (also in \cite{Giroux91}) gives a
local classification result for contact structures in a tubular
neighborhood of a convex surface.
\begin{theorem}\label{th:local_class}
If $\Sigma\neq S^2$ is a convex surface for $\xi$,
then $\Sigma$ has a tight neighborhood in $M$ if and only if no component of
$\Gamma_\Sigma$ is contractible in $\Sigma$. If $\Sigma=S^2$, then $\Sigma$
has a tight neighborhood if and only if  $\Gamma_\Sigma$ has only one
component.
\end{theorem}
Our objective is to show that, for a special choice of a Riemannian
metric, in a tubular neighborhood of $\Sigma$ the dividing set
$\Gamma_\Sigma$ becomes the set of nodal lines for a $\Delta_\Sigma$-
eigenfunction on $\Sigma$. In the next section we discuss a metric
adaptation to contact structures.

\subsection{Adapted metrics for contact structures.}
\label{sec:Riemannian metrics adapted to contact structures and their properties}
\medskip
\begin{definition}
A given metric $g$ is a {\it contact metric} for a contact form $\alpha$
if it satisfies
\begin{gather}
\label{eq:globally_adapted}
    d\alpha= \ast\lambda\,\alpha,\qquad g(\alpha,\alpha)\neq 0,\quad
\lambda\in C^{\infty}(M),\ \lambda(x)\neq 0\,\quad \text{for all }x\in M.
\end{gather}
where $\ast$ is the Hodge star operator induced by $g$.
\end{definition}
In \cite{Hamilton85}, the authors prove that any contact form $\alpha$
admits such a metric (the definitions used there are slightly
stronger).
\begin{lemma}
  Any 1-form $\alpha$ satisfying condition
\eqref{eq:globally_adapted} for some contact metric $g$ is a contact form.
\end{lemma}
\begin{proof}
One checks the contact condition
\eqref{eq:contact_condition}. We have
  \begin{gather*}
   \alpha\wedge d\alpha=\alpha\wedge(\lambda\ast\alpha)=\lambda
   g(\alpha,\alpha)\mu,\qquad \mu=\ast 1.
  \end{gather*}
By assumptions in \eqref{eq:globally_adapted} we obtain $\alpha\wedge
d\alpha\neq 0$.
\end{proof}
Every contact metric is, in fact, fully determined by a choice of an adapted
(co)frame, and can expressed in terms of a contact form $\alpha$ and its
differential $d\alpha$ \cite{Komendarczyk_thesis}.

\subsection{Nodal lines and dividing curves of contact structures.}
\qquad Recall from the introduction that the dividing set $\Gamma_\Sigma$
of a convex surface $\Sigma$ embedded in $(M,\xi)$ is the set of all points
$p$ where the contact field $v_p$ belongs to contact planes $\xi_p$.

\begin{lemma}\label{th:dividing_nodal}
Let $\Sigma$ denote a closed surface and let $\alpha$ be a contact form
on $\Sigma\times(-1,1)$ such that each $\Sigma\times\{t\}$ is convex with
a contact field $v$ preserving $\alpha$, i.e., $\Lie_v\alpha=0$.
Assume furthermore that $g$ is a contact metric satisfying
\begin{itemize}
\item[(i)] $\lambda=\text{const}$,
\item[(ii)] for each $t\in (-1,1)$, $v$ is orthonormal to
$\Sigma\times\{t\}$ with respect to $g$.
\end{itemize}
Then the dividing set $\Gamma_\Sigma$ of $\alpha$ is precisely the nodal
set of the $(-\lambda^2)$-eigenfunction of $\Delta_\Sigma$ on
$(\Sigma,g_\Sigma)$, where $g_\Sigma$ is the induced metric on $\Sigma$.
\end{lemma}
\begin{proof}
By assumption, we can choose a coframe $\{\theta_1, \theta_2, \theta_3\}$,
such that $\theta_1=g(v,\cdot)$, and $\{\theta_2, \theta_3\}\in
\Omega^1(\Sigma)$ is an orthonormal coframe on $\Sigma$. Denote by
$\{X_1,X_2,X_3\}$ a dual frame ($v=X_1$.) We can express $\alpha$ in the
coframe as follows,
  \begin{gather*}
   \alpha=f\,\theta_1 + \beta,\qquad\text{where}\quad
   \beta=\phi_2\,\theta_2+\phi_3\,\theta_3.
  \end{gather*}
  Notice that,
  \begin{gather*}
  \Gamma_\Sigma=\{p\in\Sigma\times\{0\};\ v_p\in
\xi_p\}=\{p\in\Sigma\times\{0\};\
v_p\,\lrcorner\,\alpha_p=f(p)=0\}=f^{-1}(0)\cap \Sigma\times\{0\}.
  \end{gather*}
  Now, we show that $\beta\in \Omega^1(\Sigma)$. Notice that the
requirement $\Lie_v\alpha=0$ together with \eqref{eq:globally_adapted} implies
  \begin{eqnarray*}
   0 & = &\Lie_{v}
   \alpha=v\,\lrcorner\,d\alpha+d\,f=v\,\lrcorner\,\ast\lambda\,\alpha+d\,f\\
   & = &\lambda\,v\,\lrcorner\,(
   f\,\theta_2\wedge\theta_3-\phi_2\,\theta_1\wedge\theta_3+\phi_3\,\theta_1\wedge\theta_2)+d\,f\\
   & = &\lambda\,v\,\lrcorner\,(\theta_1\wedge\ast_{\Sigma}\beta)+d\,f,
  \end{eqnarray*}
  where $\ast_\Sigma$ is the Hodge star operator on $(\Sigma,g_\Sigma)$. It
follows that,
  \begin{gather}\label{eq:lemma_nodal_1}
  d\,f=-\lambda\ast_\Sigma\,\beta.
  \end{gather}
Since $d f=\sum_i (X_i f)\, \theta_i$, we obtain the following equations
for $f$, $\phi_2$, $\phi_3$.
\begin{gather}\label{eq:coordiante_functions}
  \begin{cases}
   X_2 f=-\lambda\,\phi_3\\
   X_3 f=\lambda\,\phi_2\\
   X_1 f=0.
  \end{cases}
\end{gather}
Choosing local coordinates $(t,x,y)$, so that $v=\partial_t$ and
$(\partial_x,\partial_y)$ are tangent to the surface $\Sigma$, we conclude
that functions $f,\phi_2,\phi_3$ depend just on $(x,y)$ and
$\beta\in \Omega^1(\Sigma)$. It follows that $d \beta=d_\Sigma \beta$ and
$d\,\theta_1=d\,d\,t=0$. As a consequence of this and
\eqref{eq:lemma_nodal_1}, we obtain
$\ast\eta=\theta_1\wedge(\ast_{\Sigma} \eta)$, and
\begin{gather*}
\ast d(f\theta_1)=\ast(df\wedge\theta_1)=\lambda\ast(\theta_1\wedge\
\ast_\Sigma\beta)=\lambda\ast\ast\beta=\lambda\beta.
\end{gather*}
Expressing condition \eqref{eq:globally_adapted} in terms
of  $\alpha=f\theta_1+\beta$, we have
\begin{gather*}
  \lambda\,\alpha=\lambda\,f\theta_1+\lambda\,\beta=\ast d\alpha=\ast
d(f\theta_1)+\ast d\beta=\lambda\beta+\ast d\beta;\\
\Rightarrow\ \ast d\beta=\lambda\,f\theta_1.
\end{gather*}
It follows that
\begin{gather}\label{eq:lemma_nodal_2}
  \ast_\Sigma d\beta=\lambda\,f.
\end{gather}
Equations \eqref{eq:lemma_nodal_1} and \eqref{eq:lemma_nodal_2} imply
that $\Delta_\Sigma f=-\lambda^2\,f$,
where $\Delta_\Sigma=-\ast_\Sigma d\ast_\Sigma d$.
Therefore $\Gamma_\Sigma$ is the nodal set for $f$.
\end{proof}

Observe that in the frame $\{\theta_i\}_i$, the adapted metric $g$ is given by
\begin{gather}\label{eq:flat_metric}
  g=\sum_i \theta^2_i=dt^2+g_\Sigma,\quad \text{i.e. in coordinates}\quad
(t,x,y),\quad
  g=\begin{pmatrix}
  1 & 0\\
  0 & g_\Sigma
  \end{pmatrix},
\end{gather}
where $g_\Sigma=\theta^2_2+\theta^2_3$ is an induced metric on $\Sigma$.
Hence $g=1\oplus g_\Sigma$ is a product
metric on $U=\Sigma\times (-1,1)$. (One can prove the above lemma directly
from the decomposition of the Laplacian on $U$ in the product metric.)
\begin{remark}
If we restrict the coframe $\{\theta_i\}$ to
$\Sigma\hookrightarrow\Sigma\times\{0\}\subset \Sigma\times(-1,1)$, the
calculation in the above proof is still valid. Thus the result holds if we
assume that $\theta_1=g(v,\cdot)$ is a part of the orthonormal coframe, and
$d\,\theta_1|_{\Sigma\times\{0\}}=0$. In other words, the metric $g$ can be
only ``infinitesimally'' given as in \eqref{eq:flat_metric}.
\end{remark}
For an arbitrary $\Delta_\Sigma$-eigenfunction $f$ on $(\Sigma,
h_\Sigma)$, where $h_\Sigma$ is a smooth metric, the set of equations
\eqref{eq:coordiante_functions} determines a 1-form $\alpha$ in a
thickening of the surface $\Sigma$. If the nodal set of $f$ does not
contain singular points, then $\alpha(x)\neq 0$ for all $x$ and $\alpha$
satisfies \eqref{eq:globally_adapted} in the product metric $g=1\oplus
h_\Sigma$ (see \cite{Komendarczyk_thesis} for the explicit calculation). Consequently, we obtain the following.

\begin{theorem}
\label{th:dividing_nodal_theorem}
  If $\Xi(\Sigma,k)$ is the set of nonsingular nodal lines for a $k$-th
  eigenfunction $f$ of the Laplace operator
  $\Delta_\Sigma$ on $(\Sigma,h_\Sigma)$. Then there exists a contact
  form $\alpha$ in the thickening $\Sigma\times (-1,1)$ of  $\Sigma$ such that
  \begin{itemize}
  \item[(i)] $\Sigma$ is a convex surface for a contact structure
             $\xi_\alpha=\ker(\alpha)$,
  \item[(ii)] $\Xi(\Sigma,k)$ is the dividing set of $\xi_\alpha$.
  \item[(iii)] $\alpha$ is an $\lambda$- eigenform of the curl operator $\ast\,d$, where $-\lambda^2$ is the $\Delta_\Sigma$- eigenvalue of $f$.
  \end{itemize}
\end{theorem}
We say further that $\alpha$ is induced by $f$ in the thickening of $\Sigma$.
\subsection{The topological version of Payne's conjecture.}
 Based on results of the previous section we state the topological version of Payne's conjecture
as the following question.
\begin{problem}\label{qe:question1}
  Is the contact structure ``induced'' by the first $\Delta_\Sigma$-eigenfunctions in the sense of
  Theorem \ref{th:dividing_nodal_theorem}, always tight in the thickening of the convex surface $\Sigma$?
\end{problem}
In light of Giroux's Theorem \ref{th:local_class} the answer is positive if $\Sigma\simeq S^2$. This is a consequence of Courant's theorem, which
implies that there are exactly two nodal domains for the first
$\Delta_{S^2}$-eigenfunction on $S^2$ and that the nodal set has to be a
single embedded circle. Thus $\#\Gamma_\Sigma=\#\Xi(S^2)=1$ and the
associated contact structure has to be tight. In case $\Sigma$ is an
orientable surface of genus $\ge 1$, Problem \ref{qe:question1} is
equivalent to Question \ref{co:payne_closed} posed in the introduction. In
order to give a negative answer to \ref{qe:question1} it suffices to
construct a metric on $\Sigma$ such that the first $\Delta_\Sigma$-
eigenfunction has a closed nodal line bounding a disc. We devote the
remaining part of this paper to a rigorous construction of such metrics for
orientable surfaces of an arbitrary genus.
In \cite{Komendarczyk_thesis} we use this result, together with Theorem \ref{th:dividing_nodal_theorem}, to show existence of overtwisted principal Beltrami field (i.e. an eigenfield of the curl operator) which originally has been conjectured to be false in \cite{Ghrist00_2}.


\section{Closed nodal lines for the free membrane problem.}
\label{sec:Payne's conjecture for closed Riemannian surfaces}
\qquad Recall from the introduction that Laplace-Beltrami operator
$\Delta_M=\delta\, d$ is a positive
formally self-adjoint operator on any closed
orientable Riemannian manifold $(M,g)$. By the standard spectral theory of
formally self-adjoint
operators, the $L^2$- spectrum of $\Delta_M$ is real and countable,
\begin{gather*}
  0=\lambda_0(M)<\lambda_1(M)\leq\lambda_2(M)
  \leq\dots\leq\lambda_k(M)\leq\dots\to \infty,
\end{gather*}
and one can choose an orthonormal basis of
eigenvectors  $\{f_i\}_{i\in \mathbb{N}\cup \{0\}}$ in $L^2(M)$ (which
are smooth by regularity),
\begin{gather*}
  \Delta_M f_k=\lambda_k(M) f_k,\qquad \|f_k\|_{L^2(M)}=1,\qquad f_k\in
C^\infty(M).
\end{gather*}
The main objective of this section is to prove
\begin{theorem}\label{th:payne}
For an arbitrary closed compact orientable surface $\Sigma$, there always
exists a smooth metric $g_\Sigma$ such that $\Xi(\Sigma)$ is a single
embedded circle which bounds a disc
in $\Sigma$.
\end{theorem}
\begin{figure}
\begin{center}
  \includegraphics[width=4.5in]{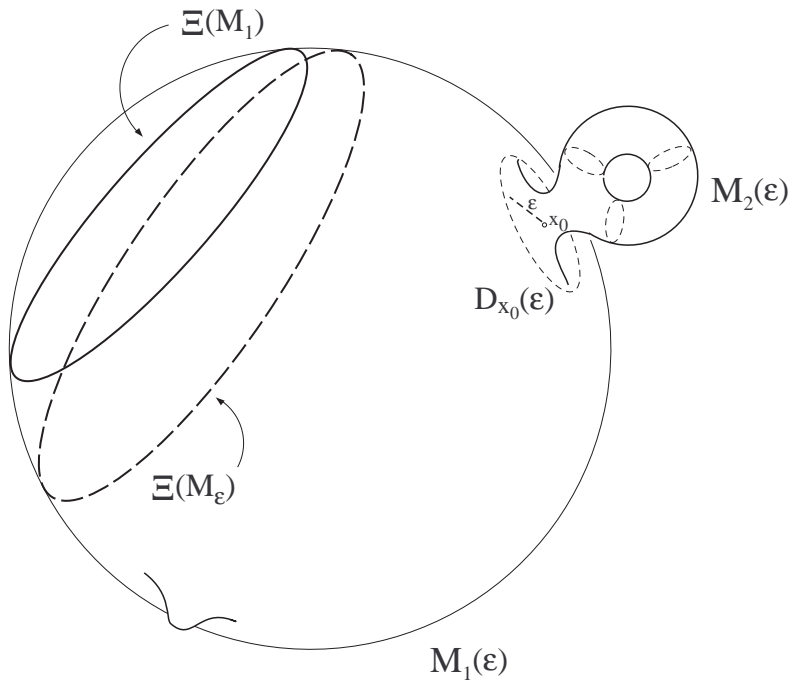}
\caption{For small $\varepsilon$, nodal lines $\Xi(M_\varepsilon)$ and
    $\Xi(M_1)$ have to be ``close'' in
    $M_\varepsilon=(M_1(\varepsilon)\cup_{\Phi_\varepsilon}
    M_2(\varepsilon),g_\varepsilon)$.}
\label{fig:nodallines}
\end{center}
\end{figure}
As already argued in the previous section, this statement is a
straightforward corollary in the case of $\Sigma=S^2$. Namely, it
is enough to choose a generic metric and refer to the result of K.
Uhlenbeck (\cite{Uhlenbeck76}) which states that $\Xi(S^2)$ has to be a one
dimensional submanifold. By Courant's Theorem $\Xi(S^2)$ splits $S^2$ into
two open domains, thereby implying that $\Xi(S^2)$ must be a single
embedded circle. If the surface is of genus $\ge 1$, we produce a desired
metric by gluing via boundary circles a ``big'' sphere $M_1\simeq S^2$ with
an $\varepsilon$- disc removed, $M_1(\varepsilon)\simeq S^2\setminus
\text{Int}(D^2)\simeq D^2$, to an $\varepsilon$-``small'' surface
$M_2(\varepsilon)$, homeomorphic to $\Sigma\setminus \text{Int}(D^2)$. The
resulting manifold
$M_\varepsilon$ is homeomorphic to $\Sigma$, and as $\varepsilon\to 0$ we
show that the nodal set $\Xi(M_\varepsilon)$ converges to $\Xi(M_1)$. (See
Figure \ref{fig:nodallines}.) Thus for sufficiently small
$\varepsilon=\hat{\varepsilon}$, $\Xi(M_{\hat{\varepsilon}})$
has to be a closed embedded circle that belongs to
$M_1(\hat{\varepsilon})\simeq D^2$.

\subsection{Definition of $M_\varepsilon$.}
\qquad First, observe the following elementary construction. If we choose
an embedded contractible 2-disc $D^2$ in an orientable surface $\Sigma$
and define $\Sigma'\cong \Sigma\setminus \text{Int}(D^2)$, then, for an
arbitrary diffeomorphism  $\Phi:\partial D^2\to\partial \Sigma'$, we can
always form a topological manifold $\Sigma_\Phi=D^2\cup_\Phi \Sigma'$ by
gluing $D^2$ back to $\Sigma'$ via $\Phi$ (see e.g. \cite{Gauld82}). Since
$\Sigma_\Phi$ is homeomorphic to $\Sigma$, we can make $\Sigma_\Phi$ into a
smooth manifold by pulling back the differential structure from $\Sigma$.
All $\Sigma_\Phi$ obtained this way are diffeomorphic. If we equip $D^2$
and $\Sigma'$ with smooth Riemannian metrics $g'_1$ and $g'_2$ we can
define a piecewise smooth metric $g$ on $\Sigma_\Phi$ as follows
  \begin{gather*}
   g=
   \begin{cases}
    g'_1 & \text{on }D^2,\\
    g'_2& \text{on }\Sigma'.
   \end{cases}
  \end{gather*}
Now $g$ is continuous on $\Sigma_\Phi$ if the gluing map $\Phi$ is an isometry.
In case $\Phi$ admits an extension to the smooth isometry of tubular
neighborhoods of boundaries $\partial D^2$, $\partial \Sigma'$, the metric
$g$ is smooth as well.

Consider an arbitrary smooth metric on $S^2$ which is flat around
$x_0$. By adding a small perturbation with support away from $x_0$, we can
produce a generic metric $g_1$ on $S^2$ flat in a small neighborhood $U_{x_0}$ of
the point $x_0$, and such that all the eigenvalues
$\{\lambda_k(M_1)\}_k$ are simple in $g_1$.
(Consult \cite{Uhlenbeck76} p. 1074 for a precise definition of a metric
perturbation and a rigorous proof of this fact in Theorem 8, p. 1076.)
 Let $M_1=(S^2,g_1)$ be a 2-sphere obtained via this
process. Additionally,   We
assume that $x_0\notin \Xi(M_1)$; otherwise we choose a different point in
the flat neighborhood. Let $D^2_{x_0}(\varepsilon)\subset U_{x_0}$ be a
geodesic disc around $x_0$ of radius $\varepsilon$ which is smaller than a
geodesic distance between $x_0$ and $\Xi(M_1)$, then for any
$\varepsilon>0$, $M_1(\varepsilon)=(M_1\setminus
\text{Int}(D^2_{x_0}(\frac{\varepsilon}{2})),g_1)$ is diffeomorphic to
$D^2$. In order to obtain a metric on $\Sigma'$, we do not make any extra
assumptions, we simply choose an arbitrary smooth metric $g_2$ on $\Sigma$,
flat around a given point $x_1$, and a geodesic disc $D^2_{x_1}(r)$ of
radius $r$ which belongs to the flat neighborhood. Clearly,
$\Sigma\setminus D^2_{x_1}(r)$ is diffeomorphic to $\Sigma'$ and  since the
metric $g_2$ can be always rescaled, we may assume that $r=1$. Define
$M_2(\varepsilon)=(\Sigma',\frac{\varepsilon^2}{4} g_2)$.

 For any $\varepsilon>0$, choose local coordinates $(x,y)$ such that the
geodesic disc $D^2_{x_0}(\varepsilon)$ is an $\varepsilon$ disc on
$(\real^2,d^2s)$ and $D^2_{x_1}(\varepsilon)$ is a unit disc on
$(\real^2,\varepsilon^2 d^2s)$, where $ds^2=dx^2+dy^2$.
Observe that the boundaries $\partial M_2(\varepsilon)$, $\partial
M_2(\varepsilon)$ can be glued via an isometry $\Phi_\varepsilon$ of
$(\real^2,d^2s)$ and $(\real^2,\frac{\varepsilon^2}{4} d^2s)$ restricted to
a circle of radius $\varepsilon$ in $(\real^2,d^2s)$. (The isometry
$\Phi_\varepsilon$ can be defined as $\Phi_\varepsilon:x\to
\frac{2}{\varepsilon} x$.)
  By the discussion in the first paragraph of this section we can form a
smooth manifold $M=M_1(\varepsilon)\cup_{\Phi_{\varepsilon}}
M_2(\varepsilon)$ and define a piecewise smooth continuous metric on $M$ as
follows (see also \cite{Takahashi02}.)
  \begin{gather}\label{eq:g_e}
   \widetilde{g}_\varepsilon = g_{M_1(\varepsilon)\cup_{\Phi_{\varepsilon}}
M_2(\varepsilon)}=
   \begin{cases}
    g_1 & \text{on }D^2,\\
   \frac{\varepsilon^2}{4} g_2& \text{on }\Sigma'.
   \end{cases}
  \end{gather}
According to \cite{Takahashi02} the following convergence of eigenvalues holds,
\begin{theorem}\label{th:taka}
For all $k=0,1,\dots$, we have
\begin{gather}
  \lim_{\varepsilon\to 0}
  \lambda_k(M,\widetilde{g}_\varepsilon)=\lambda_k(M_1,g_1).
\end{gather}
\end{theorem}
\begin{remark}
Our main objective is to prove convergence of nodal lines in $M$ to the
nodal line of $M_1$. Since the metric $\widetilde{g}_\varepsilon$ is not
smooth, and we would like to show a smooth counterexample to our version of
Payne's conjecture, we need to perturb $\widetilde{g}_\varepsilon$ in a
suitable way.
\end{remark}
For piecewise smooth metrics, eigenvalues of the Laplacian ``vary''
continuously with respect to the $C^0$-topology (see \cite{Bando83} p.
162.) Therefore, for a given $\varepsilon>0$  we can perturb the metric
$\widetilde{g}_\varepsilon$ to a smooth metric $g_\varepsilon$ so that
eigenvalues are arbitrarily ``close''. (See Theorem 1.2 in
\cite{Takahashi02}.) By Theorem 8 in \cite{Uhlenbeck76} and results of
\cite{Bando83}, we may assume that the support of the perturbation is
contained in the complement $M_1(\varepsilon)^c$ of $M_1(
\varepsilon)$ in $M$. Denote $(M,g_\varepsilon)$ by
$M_\varepsilon$.  Consequently, we can have a family of metrics
$\{g_\varepsilon\}_\varepsilon$, satisfying the following requirements.
\begin{itemize}
\item[(i)] $g_\varepsilon$ are smooth and converge to
$\widetilde{g}_\varepsilon$ in the $C^0$-topology of $M$.
\item[(ii)] $g_\varepsilon|_{M_1(\varepsilon)}=g_1$.
\item[(iii)] Eigenvalues $\lambda_k(M_\varepsilon)$ are all simple and
nodal lines $\Xi(M_\varepsilon)$ are embedded circles.
\item[(iv)]
  $\lim_{\varepsilon\to 0}\lambda_k(M_\varepsilon)=\lambda_k(M_1).$
\end{itemize}
We summarize our notation below,
\begin{itemize}
\item $M=M_1(\varepsilon)\cup_{\Phi_{\varepsilon}} M_2(\varepsilon)$,
\item $M_1=(S^2,g_1)$, $M_1(\varepsilon)=(M_1\setminus
D^2_{x_0}(\frac{\varepsilon}{2}),g_1)$,
\item $M_2(\varepsilon)=(\Sigma\setminus
D^2_{x_1}(1),\frac{\varepsilon^2}{4} g_2)$,
\item $M_\varepsilon=(M,g_\varepsilon)$.
\end{itemize}
If we must specify a different metric on a manifold, we write e.g.
$(M_2(\varepsilon),\hat{g})$.

\subsection{$C^\infty$-convergence of eigenfunctions.}
Comparing nodal lines $\Xi(M_\varepsilon)$ and $\Xi(M_1)$ can be a little
bit subtle. Notice that for each $\varepsilon>0$, $M_\varepsilon$ is
diffeomorphic to $\Sigma$ and $\{g_\varepsilon\}_\varepsilon$ is a family
of metrics on $\Sigma$. In the limit $\varepsilon=0$, the metric
$g_\varepsilon$ degenerates on $M_2(1)$, and $M_0=(\Sigma,g_0)$ is not
homeomorphic to $M_1=(S^2,g_1)$. Rather, it inherits topology that is
pulled back from $M_1$ under the quotient map, $\pi:\Sigma\to
\Sigma/M_2(1)\simeq S^2$. Thus we really have no control over what happens
to the nodal set in the ``shrinking'' part $M_2(\varepsilon)$ of the
manifold $M_\varepsilon$; technically, we cannot compare eigenfunctions on
$M_1$ to the eigenfunctions $f_k^{\varepsilon}\in C^\infty(M_\varepsilon)$
on $M_\varepsilon$. We must instead restrict them to the common domain
$M_1(\varepsilon_0)$ for a given $\varepsilon_0>0$. In order to prove the
convergence of nodal lines as $\varepsilon\to 0$, we must show  uniform
convergence of eigenfunctions $f_k^{\varepsilon}$ restricted to
$M_1(\varepsilon_0)$. In this section, we show that for any sequence
$\{\varepsilon_j\}_j$; $\varepsilon_j\to
0$,  $\{f^{\varepsilon_j}_k|_{M_1(\varepsilon_0)}\}_{\varepsilon_j}$
converges to $f_k|_{M_1(\varepsilon_0)}\in C^\infty(M_1(\varepsilon_0))$ in
the $C^\infty$-topology. For this purpose we need suitable extensions of
$f^{\varepsilon}_k|_{M_1(\varepsilon)}$ to the entire $M_1$. In the proof
we use the following extension lemma.
\begin{lemma}\label{th:main_cor}
   Let $g$ we an arbitrary smooth metric on $M_1$, given $u\in
H^l(M_1(\varepsilon),g)$, there exists a function
  $\bar{u}\in H^l(M_1,g)$, which is an extension of $u$, i.e.
  $\bar{u}|_{M_1(\varepsilon)}=u$ such that
  \begin{gather}\label{eq:extension}
  \|\bar{u}\|_{H^l(M_1,g)}\leq C\|u\|_{H^l(M_1(\varepsilon),g)}.
  \end{gather}
  For $l>1$ the constant $C$ depends on $l$ and $\varepsilon$. In the case
$l=1$ and dim=2, we can find an extension $\bar{u}$  such that $C$ is
independent of $\varepsilon$.
\end{lemma}
The proof for $l=1$ and dim=$2$ is given in \cite{Rauch75} p. 40, where the
authors show that for the unique harmonic extension, the constant $C$ is
independent of $\varepsilon$. In case $l>1$, the proof follows from the
standard extension results for functions in $H^l$ (see e.g.
\cite{Adams75}); the constant $C$ usually, grows like $1/\varepsilon^l$.)
\begin{theorem}\label{th:main}
  For each $k$, and an arbitrary $j>2$, the following $C^j$-convergence of
  eigenfunctions $f^{\varepsilon}_k\in C^{\infty}(M_\varepsilon)$,
$\|f^\varepsilon_k\|_{L^2(M_\varepsilon)}=1$ holds:
  \begin{gather}\label{eq:conv_main}
   \lim_{\varepsilon\to 0} f^{\varepsilon}_k=f_k\qquad\text{ on compact
subsets of\quad }M_1\setminus\{x_0\}
  \end{gather}
  where $f_k\in C^\infty(M_1)$ is a $k$th- $\Delta_{M_1}$-eigenfunction
on $M_1$.
\end{theorem}
  The above theorem leads to an immediate corollary,
\begin{corollary}
  Convergence \eqref{eq:conv_main} holds on compact subsets of
$M_1\setminus\{x_0\}$ in the $C^\infty$-topology of $M_1$.
\end{corollary}
\begin{proof}[Proof of Theorem \ref{th:main}]
For a given eigenfunction $f_k^{\varepsilon}$ on $M_\varepsilon$, we
introduce the following notation;
$f^\varepsilon_k=(f^{1,\varepsilon}_k,f^{2,\varepsilon}_k)$, where
$f^{1,\varepsilon}_k=f^\varepsilon_k|_{M_1(\varepsilon)}$,
$f^{2,\varepsilon}_k=f^\varepsilon_k|_{M_2(\varepsilon)}$. For convenience
we also assume $\varepsilon<1$. First we prove that there is a family of
extensions $\{\hat{f}^{1,\varepsilon}_k\}_\varepsilon$,
$\hat{f}^{1,\varepsilon}_k|_{M_1(\varepsilon)}=f^{1,\varepsilon}_k$,
convergent in $L^2(M_1)$ to $f_k$. (The argument is essentially the same as
in \cite{Takahashi02} p. 206.)

Choosing $\hat{f}^{1,\varepsilon}_k$ to be the $H^1$-extensions of
$f^{1,\varepsilon}_k\in C^\infty(M_1(\varepsilon))$ given by Lemma
\ref{th:main_cor}. We have the following,
\begin{gather}\label{eq:f_hat_1}
  \|\hat{f}^{1,\varepsilon}_k\|_{H^1(M_1)}\leq \hat{C}
\|f^{1,\varepsilon}_k\|_{H^1(M(\varepsilon),g_1)}
\end{gather}
where $\hat{C}$ is independent of $\varepsilon$. From \eqref{eq:f_hat_1} we
obtain
\begin{eqnarray}\label{eq:f_hat_2}
  \notag\|\hat{f}^{1,\varepsilon}_k\|_{H^1(M_1)}
  & \leq & \hat{C}\|f^\varepsilon_k\|_{H^{1}(M_1(\varepsilon),g_1)}
  \stackrel{(1)}{\leq} \hat{C} \bigl(\|f^\varepsilon_k\|_{L^2(M_\varepsilon)}+
                             \|d f^\varepsilon_k\|_{L^2(\Lambda^1
M_\varepsilon)}\bigr)\\
 & = &\hat{C} \bigl(1+(\Delta_{M_\varepsilon}
      f^{\varepsilon}_k,f^{\varepsilon}_k)^{\frac{1}{2}}_{L^2(M_\varepsilon)}\bigr)
\stackrel{(2)}{=}\hat{C}(1+\lambda^{\frac{1}{2}}_k(M_1)+\eta_\varepsilon),
\end{eqnarray}
where $\eta_\varepsilon\to 0$ as $\varepsilon\to 0$. Inequality (1) follows
from the definition of the $H^1$-norm, assumption (ii) on $g_\varepsilon$
and the fact that $M_1(\varepsilon)\subset M_\varepsilon$; the equality (2)
is a consequence of assumption (iv) on $g_\varepsilon$. We conclude that
the family $\{\hat{f}^{1,\varepsilon}_k\}_\varepsilon$ is bounded in
$H^1(M_1)$, thus any sequence in the family contains a weakly convergent
subsequence in $H^1(M_1)$. By Rellich's Theorem, the inclusion
$H^1(M_1)\hookrightarrow L^2(M_1)$ is compact, thus any sequence in the
family $\{\hat{f}^{1,\varepsilon}_k\}_{\varepsilon}$ contains a strongly
convergent subsequence $\{\hat{f}^{1,\varepsilon_i}_k\}_{\varepsilon_i}$ in
$L^2(M_1)$. Denote the limit of the subsequence by $\hat{f}_k\in H^1(M_1)$.
We wish to show that $\hat{f}_k$ is a smooth classical solution to
$\Delta_{M_1} u=\lambda_k(M_1) u$.

Let $B:H^1(M_1)\times H^1(M_1)\to \real$ be the bilinear form associated to
the Laplacian $\Delta_{M_1}$. Recall that for smooth functions $u$ and $w$,
\begin{gather*}
  B(u,w)=(\delta\, d\, u,w)_{L^2(M_1)}
    =(d\,u,d\,w)_{L^2(\Lambda^1 M_1)}.
\end{gather*}
The last equality extends the definition of $B$ to $H^1(M_1)$. Letting
$v\in C_c^\infty(M_1\setminus \{x_0\})$ be a test function,
\begin{eqnarray*}
B(\hat{f}_k,v)
& = &\int_{M_1} <d \hat{f}_k, d v>_{g_1} dg_1
    \stackrel{(1)}{=}\lim_{i\to\infty}\int_{M_1(\varepsilon_i)} <d
f^{1,\varepsilon_i}_k, d v>_{g_1} d g_1
\\
\\
& = &
\lim_{i\to\infty}\int_{M_1(\varepsilon_i)}
    <d f^{1,\varepsilon_i}_k, d v>_{g_1} d g_1
    +\int_{M_2(\varepsilon_i)}
    <d f^{2,\varepsilon_i}_k, 0>_{g_{\varepsilon_i}}
        d g_{\varepsilon_i}
\\
& \stackrel{(2)}{=} &
\lim_{i\to\infty} (\Delta_{M_\varepsilon}(f^{1,\varepsilon_i}_k,
        f^{2,\varepsilon_i}_k),(v,0))_{L^2(M_{\varepsilon_i})}
\\
& = &
\lim_{i\to\infty} \lambda_k(M_{\varepsilon_i})((f^{1,\varepsilon_i}_k,
        f^{2,\varepsilon_i}_k),(v,0))_{L^2(M_{\varepsilon_i})}
\\
& \stackrel{(3)}{=} &
\lambda_k(M_1)\lim_{i\to\infty}\int_{M_1(\varepsilon_i)}
    \hat{f}_k\,v\,d g_1=\lambda_k(M_1)(\hat{f}_k,v)_{L^2(M_1)}.
\end{eqnarray*}
Equality (1) follows from the $H^1$-weak convergence of extensions
$\hat{f}^{1,\varepsilon}_k$ and
$\hat{f}^{1,\varepsilon}_k|_{M_1(\varepsilon)}=f^{1,\varepsilon}_k$.
Equalities (2) and (3) follow from the assumptions (ii) and (iv) on
$g_{\varepsilon_i}$. From the density of $C^\infty_c(M_1\setminus\{x_0\})$
in $H^1(M_1,g_1)$, which holds in dimensions $\ge 2$, (see \cite{Anne86},)
the equality $B(\hat{f}_k,v)=\lambda_k(M_1)(\hat{f}_k,v)_{L^2(M_1)}$ is
valid for any $v\in H^1(M_1)$. Consequently, $\hat{f}_k$ is a weak solution
to $\Delta_{M_1} u=\lambda_k(M_1) u$, and by the regularity of weak
solutions we conclude that $\hat{f}_k$ is a smooth classical solution.
Since all the eigenvalues $\lambda_k(M_1)$ are simple in $g_1$ and
   \begin{gather*}
    1=\|f_k\|_{L^2(M_1,g_1)}=\lim_{i\to\infty}
\|f_k\|_{L^2(M_1(\varepsilon_i),g_1)}
    =\lim_{i\to\infty}
\|\hat{f}^{1,\varepsilon_i}_k\|_{L^2(M_1(\varepsilon_i),g_1)}
    =\|\hat{f}_k\|_{L^2(M_1,g_1)}
   \end{gather*}
we conclude that $\hat{f}_k=f_k$.

In the next step, we argue $C^j$-convergence of $f^{1,\varepsilon_i}_k$ on
compact subsets of $M_1\setminus\{x_0\}$.

Choose $\varepsilon_0$ such that $M_1(\varepsilon_0)$ contains a given
compact subset and
  let $l>j+\frac{m}{2}+1=j+2$, where $m=2$ is the dimension of $M_1$, (we assume $j>2$ for convenience).
Letting $\varepsilon\leq\varepsilon_0$, we apply Lemma \ref{th:main_cor}
and consider a family of $H^l$-extensions $\bar{f}^{1,\varepsilon}_k\in
H^l(M_1)$ for $f^{1,\varepsilon}_k\in C^\infty(M_1(\varepsilon_0))$. By
Lemma \ref{th:main_cor} and assumption (ii) on $g_\varepsilon$,
  \begin{gather}\label{eq:bound10}
   \|\Delta_{M_1}\bar{f}^{1,\varepsilon}_k\|_{H^{l-2}(M_1,g_1)}\leq C_{l,\varepsilon_0}\,
   \|\Delta_{M_\varepsilon} f^{1,\varepsilon}_k\|_{H^{l-2}(M_1(\varepsilon_0),g_\varepsilon)}.
  \end{gather}
Applying Garding's inequality (\cite{Roe98}, p. 76) for differential forms
with a constant $D_l$ we have
  \begin{eqnarray}
  \notag\|\bar{f}^{1,\varepsilon}_k\|_{H^l(M_1)}
  & \leq & D_l\bigl(\|\bar{f}^{1,\varepsilon}_k\|_{H^{l-1}(M_1)}
+\|(d+\delta)\bar{f}^{1,\varepsilon}_k\|_{H^{l-1}(\Lambda^\ast M_1)}\bigr)\\
\label{eq:bound1} & = & K_{l-1}\bigl(\|\bar{f}^{1,\varepsilon}_k\|_{H^{l-1}(M_1)}+\|d
  \bar{f}^{1,\varepsilon}_k\|_{H^{l-1}(\Lambda^1 M_1)}\bigr)
  \end{eqnarray}
 (where $K_{l-1}=D_l$.) Here $d+\delta$ is the Dirac operator (i.e. $(d+\delta)^2=\Delta$) acting
on forms of mixed degree. Applying Garding's inequality again to each term
of \eqref{eq:bound1} and setting $D_{l,l-1}=D_l D_{l-1}$ results in
  \begin{eqnarray*}
  \text{(rhs of \eqref{eq:bound1})}
  &\leq &
  D_{l,l-1} \bigl(\|\bar{f}^{1,\varepsilon}_k\|_{H^{l-2}(M_1)}+
  2 \|d\bar{f}^{1,\varepsilon}_k\|_{H^{l-2}(\Lambda^1 M_1)}+
  \|\Delta_{M_1} \bar{f}^{1,\varepsilon}_k\|_{H^{l-2}(M_1)}\bigr)
  \\
  & \stackrel{(1)}{\leq} &
   D_{l,l-1} \bigl(\|\bar{f}^{1,\varepsilon}_k\|_{H^{l-2}(M_1)}+
   2 \|d\bar{f}^{1,\varepsilon}_k\|_{H^{l-2}(\Lambda^1 M_1)}+
   C_{l-2,\varepsilon_0}\|\Delta_{M_\varepsilon} f^{\varepsilon}_k\|_{H^{l-2}(M_1(\varepsilon_0),g_\varepsilon)}\bigr)
  \\
  &\stackrel{(2)}{\leq}&
  D_{l,l-1}
  \bigl((1+\lambda_k(M_\varepsilon)C_{l-2,\varepsilon_0})\|\bar{f}^{1,\varepsilon}_k\|_{H^{l-2}(M_1)}+
   2 \|d \bar{f}^{1,\varepsilon}_k\|_{H^{l-2}(\Lambda^1 M_1)}\bigr)\\
  &\stackrel{(3)}{\leq}&
K_{l-2}\bigl(\|\bar{f}^{1,\varepsilon}_k\|_{H^{l-2}(M_1)}+
   \|d\bar{f}^{1,\varepsilon}_k\|_{H^{l-2}(\Lambda^1 M_1)}\bigr)
  \end{eqnarray*}
 Inequality (1) is a consequence of \eqref{eq:bound10}, whereas (2)
follows from the fact that
$f^{\varepsilon}_k$ is the $k$th eigenfunction of the Laplacian on
$M_\varepsilon$. In (3) we can set $K_{l-2}> D_{l,l-1} \max\{2,(1+\lambda_k(M_\varepsilon)C_{l-2, \varepsilon_0})\}$ due to the requirement (iv) on $g_\varepsilon$.
 Repeating the above steps finitely many times leads to
   \begin{eqnarray*}
     \|\bar{f}^{1,\varepsilon}_k\|_{H^l(M_1)}& \leq &
K_1  \|\bar{f}^{1,\varepsilon}_k\|_{H^1(M_1)}\leq
K_1  C_1\bigl(\|f^{1,\varepsilon}_k\|_{L^2(M_1(\varepsilon_0),g_1)}+
     \|df^{1,\varepsilon}_k\|_{L^2(\Lambda^1 M_1(\varepsilon_0),g_1)}\bigr)\\
   & \leq & K_1 C_1\bigl(1+(\Delta_{M_\varepsilon}
f^{\varepsilon}_k,f^{\varepsilon}_k)^{\frac{1}{2}}_{L^2(M_\varepsilon)}\bigr)
   \leq K_1 C_1 (1+\lambda^{\frac{1}{2}}_k(M_1)+\eta'_\varepsilon),
   \end{eqnarray*}
   where $\eta'_\varepsilon\to 0$ as $\varepsilon\to 0$ and we applied
Lemma \ref{th:main_cor} in the second inequality. Consequently, the family
$\{\bar{f}^{1,\varepsilon}_k\}_\varepsilon$ is bounded in ${H^l(M_1)}$. By
Rellich's Theorem we have a compact inclusion $H^l(M_1)\hookrightarrow
H^{l-1}(M_1)$, and by Sobolev embedding theorem a bounded inclusion
$H^{l-1}(M_1)\hookrightarrow C^{j}(M_1)$. Composition of these two gives us
a compact inclusion  $H^l(M_1)\hookrightarrow C^{j}(M_1)$. As a result,
there exists a subsequence $\{\bar{f}^{1,\varepsilon_i}_k\}_i$ of any
sequence in $\{\bar{f}^{1,\varepsilon}_k\}_\varepsilon$, convergent in the
$C^{j}$ topology of $M_1$. Denote a limit of this  subsequence by
$\bar{f}_k\in C^j(M_1)$.

Since
$\hat{f}^{1,\varepsilon}_k|_{M_1(\varepsilon_0)}=f^{1,\varepsilon}_k|_{M_1(\varepsilon_0)}=\bar{f}^{1,\varepsilon}_k|_{M_1(\varepsilon_0)}$
for any $\varepsilon<\varepsilon_0$, the $L^2$-limits $\hat{f}_k$,
$\bar{f}_k$ have to agree on $M_1(\varepsilon_0)$. Equality
\eqref{eq:conv_main} holds, since for any sequence $\{\varepsilon_i\}_i$
converging to zero $\{f_k^{1,\varepsilon_i}|_{M_1(\varepsilon_0)}\}$
contains a convergent subsequence with a common limit.
\end{proof}

\subsection{Proof of Theorem \ref{th:payne}.}\label{sec:Proof of Theorem}
We assume that the nodal line $\Xi(M_1)$ is at a geodesic distance
$d>\varepsilon_0$ from the gluing disc $D_{\varepsilon}\subset
D_{\varepsilon_0}$. The nodal sets $\Xi(M_1)$ and $\Xi(M_\varepsilon)$ can
be compared only on the common subset $M_1(\varepsilon_0)$, and since we
just proved the convergence in $C^\infty(M_1(\varepsilon_0))$ of
$f^{\varepsilon_i}_1|_{M_1(\varepsilon_0)}\to f_1|_{M_1(\varepsilon_0)}$
for any sequence $\varepsilon_i\to 0$, the first impression is that we have
no control over nodal lines in the ``shrinking'' part $M_2(\varepsilon)$ of
the manifold $M_\varepsilon$. What we really have to show is that for some
$\varepsilon>0$ the nodal set $\Xi(M_\varepsilon)$ belongs entirely to the
common domain $M_1(\varepsilon_0)$. This would imply that
$\Xi(M_\varepsilon)$ is an embedded contractible circle,
since  $M_1(\varepsilon_0)$ is itself contractible. First of all note the
following pointwise convergence of nodal sets.
\begin{lemma}\label{th:pointwise}
  Consider a sequence of points $\{x_i\}_i$ such that for each $i$,
  $x_i\in \Xi(M_{\varepsilon_i})\cap M_1(\varepsilon_0)$.
   If the limit $x$ of $\{x_i\}_i$ exists, then $x\in \Xi(M_1)$.
\end{lemma}
\begin{proof}
   Applying the convergence
   $f^{\varepsilon_i}_1|_{M_1(\varepsilon_0)}\to f_1|_{M_1(\varepsilon_0)}$
in $C^0(M_1(\varepsilon_0))$ we obtain,
  \begin{gather*}
  |f_1(x_i)|=|f_1(x_i)-f^{1,\varepsilon_i}_1(x_i)|\leq
  \|f_1-f^{1,\varepsilon_i}_1\|_{C^{0}(M_1(\varepsilon_0))}
  \stackrel{i\to \infty}{\longrightarrow} 0.
  \end{gather*}
   By continuity of $f_1$ and the assumption $x_i\to x\in
M_1(\varepsilon_0)$, we have $0=\lim_{i\to \infty}
   f_1(x_i)=f_1(x)$, and we conclude that $x\in \Xi(M_1)$.
\end{proof}
  In order to finish the proof of Theorem \ref{th:payne} we need only show
the following.
  \begin{claim}
   The nodal lines $\Xi(M_{\varepsilon_i})$ are eventually in
  $M_1(\varepsilon_0)$, i.e. we can find an index $n$ such that for all
  $i>n$, $\Xi(M_{\varepsilon_i})\subset M_1(\varepsilon_0)$.
  Additionally for $i>n$, each $\Xi(M_{\varepsilon_i})$ is a single
embedded circle.
  \end{claim}
   We consider two cases; either there exists
   an infinite sequence of points $\{x_i\}_i\subset
   \Xi(M_{\varepsilon_i})\cap M_1(\varepsilon_0)$ (case $1^\circ$), or not
(case $2^\circ$).
\begin{proof}[Proof of the claim in case $\it 1^\circ$.]
  By compactness of $M_1(\varepsilon_0)$, we can choose a convergent
subsequence of $\{x_i\}_i$. We denote the subsequence by $\{\hat{x}_j\}_j$
and its limit by $\hat{x}$. According to Lemma \ref{th:pointwise}, $\hat{x}\in
  \Xi(M_1)\subset M_1(\varepsilon_0)$. By assumption (i) on p. 10, for each
$j$ we have an embedding $\phi_j:S^1\hookrightarrow M_{\varepsilon_j}$ such
that for some $\theta_j\in S^2$, $\phi_j(\theta_j)=\hat{x}_j$. If
$\varepsilon_j$ is sufficiently small, all elements of $\hat{x}_j$ belong
to a geodesic ball $D_{\hat{x}}(\hat{r})\subset M_1(\varepsilon_0)$ around
$\hat{x}$ of arbitrarily small radius $\hat{r}$. To prove the claim, we
reason by contradiction. Suppose that there exists a subsequence $\{m\}$ of
$\{j\}$ such that $\phi_{m}(S^1)\nsubseteq M_1(\varepsilon_0)$. Then for
each $m$ there is a point $\vartheta_{m}$ where
$\phi_{m}(\vartheta_{m})\notin M_1(\varepsilon_0)$. Consequently,
$\phi_m:(\theta_m,\vartheta_m)\hookrightarrow M_{\varepsilon_m}$ is an
oriented path joining two different points, one in $M_1(\varepsilon_0)$,
the other in $M_{\varepsilon_m}\setminus M_1(\varepsilon_0)$. Clearly, the
path $\phi_m|_{(\theta_m,\vartheta_m)}$ has to intersect $\partial
M_1(\varepsilon_0)$. Choose a convergent subsequence of intersection points
$y_m\in \partial M_1(\varepsilon_0)$. Its limit $y\in \partial
M_1(\varepsilon_0)$ belongs to $\Xi(M_1)$ by Lemma \ref{th:pointwise}. As a
consequence, the intersection $\Xi(M_1)\cap \partial M_1(\varepsilon_0)$
would have to be nonempty which contradicts the choice of $\varepsilon_0$
(i.e. the nodal line $\Xi(M_1)$ is supposed to be at some geodesic distance
$d>\varepsilon_0$ from the boundary of the gluing disc.)
  Therefore, there exists an index $n$ such that $\phi_j(S^1)\subset
M_1(\varepsilon_0)$ for all $j>n$. The image of $S^1$ under $\phi_j$ is a
closed embedded curve in $M_1(\varepsilon_0)\simeq D^2$ which splits
$M_{\varepsilon_j}$ into two domains, where the eigenfunction is either
strictly positive or negative (a $\Delta_\Sigma$-eigenfunction cannot have
a zero so there is no change of sign in the neighborhood). By Courant's
Theorem, the first eigenfunction must have exactly two nodal domains, hence
$\phi_j(S^1)=\Xi(M_{\varepsilon_j})$.
\end{proof}

\begin{proof}[Proof of the claim in case $\it 2^\circ$.]
If there is no infinite sequence $\{x_i\}_i\subset\Xi(M_{\varepsilon_i})\cap
M_1(\varepsilon_0)$, then there exists an index $n$ such that for any
$i>n$, $\Xi(M_{\varepsilon_i})\subset M_2(1)$. Again by Courant's Theorem,
for each $i>n$, $M_{\varepsilon_i}\setminus \Xi(M_{\varepsilon_i})$
consists of two connected disjoint open subsets $M^{+}_{\varepsilon_i}$,
$M^{-}_{\varepsilon_i}$ of $M_{\varepsilon_i}$, defined as
$M^{+(-)}_{\varepsilon_i}=\{x\in M_{\varepsilon_i};
f^{1,\varepsilon_i}_1(x)>(<)0\}$. As a result $M_1(\varepsilon_0)\subset
M^{+}_{\varepsilon_i}$, or $M_1(\varepsilon_0)\subset
M^{-}_{\varepsilon_i}$, for all $i>n$. The convergence
$f^{1,\varepsilon_i}_1\longrightarrow f_1$ in $C^0(M_1(\varepsilon_0))$
implies that $f_1$ has to be either nonnegative on $M_1(\varepsilon_0)$ or
nonpositive. But this cannot happen, because $\Xi(M_1)\subset
M_1(\varepsilon_0)$ and $f_1$ has to change sign on $M_1(\varepsilon_0)$.
\end{proof}

\section{Conclusions}
Lemma \ref{th:dividing_nodal} ``ties'' the nodal set of an eigenfunction to
dividing curves
of an induced contact structure. It would be very interesting to show that
dividing curves of an arbitrary contact structure are nodal lines in some
suitably adapted metric.
\begin{problem}\label{pr:problem1}
Let $(M,\xi)$ be a contact 3-manifold and $\Sigma$ an embedded convex
surface. Can one adapt a metric $g$ to $\xi$ in a tubular neighborhood of
$\Sigma$ so that the dividing curves correspond to nodal lines of a
$\Delta_\Sigma$-eigenfunction on $(\Sigma,g_\Sigma)$ ?
\end{problem}
Since isotopic contact structures have isotopic dividing curves, in the
formulation of the above problem we can allow the distribution $\xi$ to
vary up to isotopy. If the answer to Problem \ref{pr:problem1} is ``yes''
one could rephrase questions about nodal lines in terms of dividing curves
and vice versa. Especially, one could try to address questions of the
following type using topological techniques.
\begin{problem}
  What kind of nonsingular nodal lines are admissible by a given compact
orientable surface $\Sigma$? In other words, can we realize a given
collections of embedded closed curves $\{\Gamma_i\}_i$ on $\Sigma$ as nodal
lines (possibly up to isotopy) for some Riemannian metric $g_\Sigma$ on
$\Sigma$.
\end{problem}
There are some obvious restrictions on the family $\{\Gamma_i\}_i$, e.g. it
must divide the surface $\Sigma$.
To generate examples of such families, we could adapt the
technique from Section \ref{sec:Payne's conjecture for  closed Riemannian
surfaces}. It seems possible to start with a generic nodal set embedded
into a 2-sphere $S^2$ and ``implant'' a finite number of small handles (as
pictured on Figure \ref{fig:nodallines}), which collapse to centers of
attaching circles, as a parameter $\varepsilon\to 0$. By similar reasoning
as in Section \ref{sec:Proof of Theorem} we could find a ``generic'' metric
on a surface obtained during this process, with the nodal set isotopic to
the initial nodal set on $S^2$. One could also start with an arbitrary
orientable surface and carry out the implanting procedure. It may also be
interesting to relate this to gluing of contact structures along convex
surfaces with Legendrian boundary (see \cite{Kazez02}).

In the context of Giroux's Theorem \ref{th:local_class}, and the active
search for the classification result for tight contact structures, it may
be interesting to address the following problem.
\begin{problem}
  Find conditions, if any, on the curvature of a closed surface so that all
the nodal lines of a $\Delta_\Sigma$-eigenfunction are homotopically
nontrivial.
\end{problem}
The hyperbolic case would seem to be of special interest here.

\small
\section{Acknowledgments}
Advice and comments from M. Dillon, R. Ghrist, J. Landsberg, M. Symington, and
A. Swiech were greatly helpful. This paper is part of the author's
Ph.D. thesis work under the supervision of R. Ghrist.

\end{document}